\documentclass{amsart}
\usepackage{amssymb}
\theoremstyle{plain}
 \newtheorem{thm}{Theorem}
 
 \newtheorem{lem}{Lemma}
 \newtheorem{prop}{Proposition}
\theoremstyle{definition}

 \newtheorem{exmp}{Example}
\theoremstyle{remark}
 \newtheorem{rem}{Remark\ignorespaces}
 
\newcommand{\NaturalNumber}{\mathbb N}
\newcommand{\RealNumber}{\mathbb R}
\renewcommand{\labelenumi}{(\roman{enumi})}
\newcommand{\toMUGEN}[3]{
 \mathop{#1_{#2 \in #3,}}_{\| #2 \| \rightarrow \infty}}

\begin{document}

\title[Asymptotically contractive mapping]
{Fixed point theorems for asymptotically contractive mappings}
\author[T. Suzuki]{Tomonari Suzuki}
\date{}
\hyphenation{kyu-shu kita-kyu-shu to-bata-ku sen-sui-cho}
\address{
Department of Mathematics,
Kyushu Institute of Technology,
1-1, Sensuicho, Tobataku, Kitakyushu 804-8550, Japan}
\email{suzuki-t@mns.kyutech.ac.jp}

\begin{abstract}
In this short paper,
 we prove fixed point theorems for nonexpansive mappings
 whose domains are unbounded subsets of Banach spaces.
These theorems are generalizations of Penot's result
 in [Proc.\ Amer.\ Math.\ Soc., 131 (2003), 2371--2377].
\end{abstract}
\maketitle

\section{Introduction}
\label{SC:intro}

Let $C$ be a closed convex subset of a Banach space $E$, and
 let $T$ be a {\it nonexpansive mapping} on $C$, i.e.,
 $\| Tx - Ty \| \leq \| x - y \|$
 for all $x, y \in C$.
We know that $T$ has a fixed point
 in the case that $E$ is uniformly convex and $C$ is bounded;
 see Browder \cite{REF:Browder1965_ProcNAS_3} and
 G\"ohde \cite{REF:Gohde1965}.
Kirk \cite{REF:Kirk1965_AMMonth} extended these result
 to the case that $C$ is weakly compact and has normal structure.
We note that
 such domain $C$ of $T$ is a bounded subset.
Recently,
 Penot proved the following in \cite{REF:Penot2003_ProcAMS}:
 $T$ has a fixed point
 in the case that $E$ is uniformly convex,
 $C$ is unbounded, and
 $T$ is {\it asymptotically contractive}, i.e.,
 $$ \toMUGEN{\limsup}{y}{C}
 \frac{\| T x_0 - T y \|}{\| x_0 - y \|} < 1 $$
 for some $x_0 \in C$.

In this paper,
 we prove fixed point theorems for nonexpansive mappings
 whose domains are unbounded subsets of Banach spaces.
These theorems are generalizations of Penot's result in
 \cite{REF:Penot2003_ProcAMS}.

\section{Conditions for Mappings}
\label{SC:condi}

In this section,
 let $T$ be a nonexpansive mapping on a nonempty closed convex subset $C$
 of a Banach space $E$.
We discuss the following conditions for $T$:
\begin{enumerate}
\renewcommand{\labelenumi}{(C\arabic{enumi})}
\renewcommand{\theenumi}{(C\arabic{enumi})}
\item\label{ENUM:CONDI:asympt-cont-a}
 There exists $r \in (0,1)$ such that
 for every $x_1 \in C$,
 there exists $\eta > 0$ satisfying
 $$ \| T x_1 - T y \| \leq r \; \| x_1 - y \| $$
 for all $y \in C$ with $\| y \| > \eta$;
\item\label{ENUM:CONDI:asympt-cont-e}
 there exist $r \in (0,1)$, $x_0 \in C$ and $\eta > 0$
 such that
 $$ \| T x_0 - T y \| \leq r \; \| x_0 - y \| $$
 for all $y \in C$ with $\| y \| > \eta$;
\item\label{ENUM:CONDI:infinity-a}
 for each $\lambda > 0$ and for each $x_1 \in C$,
 there exists $\eta > 0$ satisfying
 $$ \| T x_1 - T y \| \leq \| x_1 - y \| - \lambda $$
 for all $y \in C$ with $\| y \| > \eta$;
\item\label{ENUM:CONDI:infinity-e}
 there exists $x_0 \in C$
 for each $\lambda > 0$,
 there exists $\eta > 0$ satisfying
 $$ \| T x_0 - T y \| \leq \| x_0 - y \| - \lambda $$
 for all $y \in C$ with $\| y \| > \eta$;
\item\label{ENUM:CONDI:for-afpp}
 there exists $\lambda > 0$ such that
 for each $x_1 \in C$,
 there exists $\eta > 0$ satisfying
 $$ \| T x_1 - T y \| \leq \| x_1 - y \| - \lambda $$
 for all $y \in C$ with $\| y \| > \eta$;
\item\label{ENUM:CONDI:basic}
 there exist $x_0 \in C$ and $\eta > 0$ such that
 $$ \| T x_0 - T y \| \leq \| x_0 - y \| - \| T x_0 - x_0 \| $$
 for all $y \in C$ with $\| y \| > \eta$.
\end{enumerate}

We obtain the following.

\begin{prop}
\label{PROP:conditions}
\ref{ENUM:CONDI:asympt-cont-a}
 $\Leftrightarrow$ \ref{ENUM:CONDI:asympt-cont-e}
 $\Rightarrow$ \ref{ENUM:CONDI:infinity-a}
 $\Leftrightarrow$ \ref{ENUM:CONDI:infinity-e}
 $\Rightarrow$ \ref{ENUM:CONDI:for-afpp}
 $\Rightarrow$ \ref{ENUM:CONDI:basic}
 holds.
\end{prop}

\begin{proof}
It is obvious that
 \ref{ENUM:CONDI:asympt-cont-a}
 $\Rightarrow$ \ref{ENUM:CONDI:asympt-cont-e},
 \ref{ENUM:CONDI:infinity-a}
 $\Rightarrow$ \ref{ENUM:CONDI:infinity-e}, and
 \ref{ENUM:CONDI:infinity-a}
 $\Rightarrow$ \ref{ENUM:CONDI:for-afpp}.
We first prove
 \ref{ENUM:CONDI:asympt-cont-e}
 $\Rightarrow$ \ref{ENUM:CONDI:asympt-cont-a}.
We assume \ref{ENUM:CONDI:asympt-cont-e}, i.e.,
 there exists $r' \in (0,1)$, $x_0 \in C$ and $\eta' > 0$
 such that
 $ \| T x_0 - T y \| \leq r' \; \| x_0 - y \| $
 for all $y \in C$ with $\| y \| > \eta'$.
Put $r = (1+r')/2$.
We let $x_1 \in C$ be fixed and
 put
 $$ \eta = \max\left\{ \eta',
 \| x_1 \| + \frac{\| x_1 - x_0 \| + \| T x_1 - T x_0 \|}{r-r'} \right\} . $$
Then for $y \in C$ with $\| y \| > \eta$,
 we have
 \begin{align*}
 \| x_1 - x_0 \| + \| T x_1 - T x_0 \|
 &\leq (r - r') \; (\eta - \| x_1 \|) \\*
 &\leq (r - r') \; (\| y \| - \| x_1 \|) \\*
 &\leq (r - r') \; \| x_1 - y \|
 \end{align*}
 and hence
 \begin{align*}
 \| T x_1 - T y \|
 &\leq \| T x_1 - T x_0 \| + \| T x_0 - T y \| \\*
 &\leq \| T x_1 - T x_0 \| + r' \| x_0 - y \| \\
 &\leq \| T x_1 - T x_0 \| + r' \| x_1 - x_0 \| + r' \| x_1 - y \| \\
 &\leq \| T x_1 - T x_0 \| + \| x_1 - x_0 \| + r' \| x_1 - y \| \\*
 &\leq r \| x_1 - y \| .
 \end{align*}
This implies \ref{ENUM:CONDI:asympt-cont-a}.
We can similarly prove
 \ref{ENUM:CONDI:asympt-cont-e}
 $\Rightarrow$ \ref{ENUM:CONDI:infinity-e} and
 \ref{ENUM:CONDI:infinity-e}
 $\Rightarrow$ \ref{ENUM:CONDI:infinity-a}.
We finally show
 \ref{ENUM:CONDI:for-afpp} $\Rightarrow$ \ref{ENUM:CONDI:basic}.
We assume \ref{ENUM:CONDI:for-afpp}, i.e.,
 there exists $\lambda > 0$ such that
 for each $x_1 \in C$,
 there exists $\eta > 0$ satisfying
 $ \| T x_1 - T y \| \leq \| x_1 - y \| - \lambda $
 for all $y \in C$ with $\| y \| > \eta$.
We put $$ d = \inf_{x \in C} \| T x - x \| $$
 and assume $d > 0$.
Then there exists $x_1 \in C$ such that $\| T x_1 - x_1 \| < d + \lambda / 2$.
For such $x_1$,
 we choose $\eta > 0$ satisfying
 $ \| T x_1 - T y \| \leq \| x_1 - y \| - \lambda $
 for all $y \in C$ with $\| y \| > \eta$.
For each $t \in (0,1)$,
 since a mapping $x \mapsto (1-t) T x + t x_1$ on $C$ is contractive,
 there exists $y_t \in C$
 such that
 $$ y_t = (1 - t) T y_t + t x_1 .$$
Since
 \begin{align*}
 d
 &\leq \| T y_t - y_t \|
 = t \| T y_t - x_1 \| \\*
 &\leq t \; \big( \| T y_t - T x_1 \| + \| T x_1 - x_1 \| \big) \\
 &\leq t \; \big( \| y_t - x_1 \| + \| T x_1 - x_1 \| \big) \\
 &\leq t \; \big( \| y_t \| + \| x_1 \| + \| T x_1 - x_1 \| \big) , \\*
 \end{align*}
 we have
 $\| y_t \| > \eta$ for some small $t > 0$.
So, we have
 \begin{align*}
 \| x_1 - y_t \| + \| y_t - T y_t \|
 &= \| x_1 - T y_t \| \\*
 &\leq \| x_1 - T x_1 \| + \| T x_1 - T y_t \| \\
 &\leq \| x_1 - T x_1 \| + \| x_1 - y_t \| - \lambda \\*
 &\leq d + \lambda / 2 + \| x_1 - y_t \| - \lambda
 \end{align*}
 and hence
 $$ \| y_t - T y_t \| \leq d - \lambda / 2 .$$
This contradicts to the definition of $d$.
Therefore we obtain $d = 0$.
We can choose $x_0 \in C$ with $\| T x_0 - x_0 \| < \lambda$.
Then there exists $\eta > 0$ such that
 \begin{align*}
 \| T x_0 - T y \|
 &\leq \| x_0 - y \| - \eta \\
 &< \| x_0 - y \| - \| T x_0 - x_0 \|
 \end{align*}
 for all $y \in C$ with $\| y \| > \eta$.
This completes the proof.
\end{proof}

We can easily prove the following.

\begin{prop}
\label{PROP:unbounded}
Suppose that $C$ is unbounded.
Then the following are equivalent to \ref{ENUM:CONDI:asympt-cont-a}
 and \ref{ENUM:CONDI:asympt-cont-e}:
\begin{enumerate}
\item $T$ is asymptotically contractive;
\item for every $x_1 \in C$,
 $$ \toMUGEN{\limsup}{y}{C}
 \frac{\| T x_1 - T y \|}{\| x_1 - y \|} < 1 $$
 holds.
\end{enumerate}
And the following are equivalent to \ref{ENUM:CONDI:infinity-a}
 and \ref{ENUM:CONDI:infinity-e}:
\begin{enumerate}
\item there exists $x_0 \in C$ such that
 $$ \toMUGEN{\lim}{y}{C}
 \Big( \| T x_0 - T y \| - \| x_0 - y \| \Big) = - \infty . $$
\item for every $x_1 \in C$,
 $$ \toMUGEN{\lim}{y}{C}
 \Big( \| T x_1 - T y \| - \| x_1 - y \| \Big) = - \infty $$
 holds.
\end{enumerate}
\end{prop}

\section{Sufficient and Necessary Condition}
\label{SC:suf-nec}

In this section,
 we discuss about the sufficient and necessary condition
 for nonexpansive mappings having a fixed point.

\begin{lem}
\label{LEM:invariant}
Let $C$ be a closed convex subset of a Banach space $E$ and
 let $T$ be a nonexpansive mapping on $C$.
Suppose that \ref{ENUM:CONDI:basic}, i.e.,
 there exist $x_0 \in C$ and $\eta > 0$ such that
 $$ \| T x_0 - T y \| \leq \| x_0 - y \| - \| T x_0 - x_0 \| $$
 for all $y \in C$ with $\| y \| > \eta$.
Then there exists $\rho > 0$ such that
 $ T(D) \subset D $,
 where
 $$ D = \{ y \in C : \| y - x_0 \| \leq \rho \} .$$
\end{lem}

\begin{proof}
We put
 $$ \rho = \eta + \| x_0 \| + \| T x_0 - x_0 \| > 0 . $$
Then in the case of $y \in D$ and $\| y \| \leq \eta$,
 we have
 \begin{align*}
 \| T y - x_0 \|
 &\leq \| T y - T x_0 \| + \| T x_0 - x_0 \| \\*
 &\leq \| y - x_0 \| + \| T x_0 - x_0 \| \\
 &\leq \| y \| + \| x_0 \| + \| T x_0 - x_0 \| \\
 &\leq \eta + \| x_0 \| + \| T x_0 - x_0 \| \\*
 &= \rho .
 \end{align*}
In the case of $y \in D$ and $\| y \| > \eta$,
 we have
 \begin{align*}
 \| T y - x_0 \|
 &\leq \| T y - T x_0 \| + \| T x_0 - x_0 \| \\*
 &\leq \| y - x_0 \| \\*
 &\leq \rho .
 \end{align*}
Therefore we obtain the desired result.
\end{proof}

A closed convex subset $C$ of a Banach space $E$ is said to have
 the {\it fixed point property} for nonexpansive mappings
 ({\it FPP}, for short)
 if for every bounded closed convex subset $D$ of $C$,
 every nonexpansive mapping on $D$ has a fixed point.
Similarly,
 $C$ is said to have
 the {\it weak fixed point property} for nonexpansive mappings
 ({\it WFPP}, for short)
 if for every weakly compact convex subset $D$ of $C$,
 every nonexpansive mapping on $D$ has a fixed point.
Let $E^\ast$ be the dual of $E$.
Then
 a closed convex subset $C$ of $E^\ast$ is said to have
 the {\it weak$^\ast$ fixed point property} (with respect to $E$)
 for nonexpansive mappings
 ({\it W$^\ast$FPP}, for short)
 if for every weakly$^\ast$ compact convex subset $D$ of $C$,
 every nonexpansive mapping on $D$ has a fixed point.
So, by the results of Browder \cite{REF:Browder1965_ProcNAS_3}
 and G\"ohde \cite{REF:Gohde1965},
 every uniformly convex Banach space has FPP.
Also, by Kirk's result \cite{REF:Kirk1965_AMMonth},
 every Banach space with normal structure has WFPP.
We recall that a closed convex subset $C$ of a Banach space $E$
 is {\it locally weakly compact}
 if and only if
 every bounded closed convex subset of $C$ is weakly compact.
So, every closed convex subset of a reflexive Banach space
 is locally weakly compact.

Using Lemma \ref{LEM:invariant}, we obtain the following propositions.

\begin{prop}
\label{PROP:fpp}
Let $C$ be a closed convex subset of a Banach space $E$.
Assume that
 $C$ has FPP.
Let $T$ be a nonexpansive mapping on $C$.
Then the following are equivalent:
\begin{enumerate}
\item $T$ has a fixed point in $C$;
\item $T$ satisfies \ref{ENUM:CONDI:basic}.
\end{enumerate}
\end{prop}

\begin{proof}
We first show (ii) implies (i).
We suppose that (ii), i.e.,
 there exist $x_0 \in C$ and $\eta > 0$ such that
 $$ \| T x_0 - T y \| \leq \| x_0 - y \| - \| T x_0 - x_0 \| $$
 for all $y \in C$ with $\| y \| > \eta$.
By Lemma \ref{LEM:invariant},
 there exists $\rho > 0$ such that
 $ T(D) \subset D $,
 where
 $$ D = \{ x \in C : \| x - x_0 \| \leq \rho \} .$$
So, by the assumption,
 there exists $z_0 \in D$ such that $T z_0 = z_0$.
Conversely, let us prove that (i) implies (ii).
Let $x_0$ be a fixed point of $T$.
Since $T$ is nonexpansive,
 we have
 $$ \| T x_0 - T y \| \leq \| x_0 - y \|
 = \| x_0 - y \| - \| T x_0 - x_0 \| $$
 for all $y \in C$.
This implies \ref{ENUM:CONDI:basic}.
This completes the proof.
\end{proof}

\begin{prop}
\label{PROP:wfpp}
Let $C$ be a closed convex subset of a Banach space $E$.
Assume that
 $C$ is locally weakly compact and has WFPP.
Let $T$ be a nonexpansive mapping on $C$.
Then $T$ has a fixed point in $C$
 if and only if
 $T$ satisfies \ref{ENUM:CONDI:basic}.
\end{prop}

\begin{prop}
\label{PROP:wfpp-ast}
Let $E$ be a Banach space and
 let $E^\ast$ be the dual of $E$.
Let $C$ be a weakly$^\ast$ closed convex subset of $E^\ast$.
Assume that
 $C$ has W$^\ast$FPP.
Let $T$ be a nonexpansive mapping on $C$.
Then $T$ has a fixed point in $C$
 if and only if
 $T$ satisfies \ref{ENUM:CONDI:basic}.
\end{prop}

As a direct consequence, we have the following.

\begin{thm}
\label{THM:basic}
Let $E$ be a Banach space and let $E^\ast$ be the dual of $E$.
Assume that either of the following:
\begin{enumerate}
\item $C$ is closed convex subset of $E$ and has FPP.
\item $C$ is closed convex subset of $E$ and locally weakly compact,
 and $C$ has WFPP.
\item $C$ is weakly$^\ast$ closed convex subset of $E^\ast$
 and has W$^\ast$FPP.
\end{enumerate}
Let $T$ be a nonexpansive mapping on $C$.
Suppose that $C$ is unbounded, and $T$ is asymptotically contractive.
Then $T$ has a fixed point.
\end{thm}

\begin{rem}
(ii) implies (i).
\end{rem}

\begin{thm}[Penot \cite{REF:Penot2003_ProcAMS}]
\label{THM:Penot}
Let $C$ be a unbounded closed convex subset of
 a uniformly convex Banach space $E$.
Let $T$ be a nonexpansive mapping on $C$.
Suppose that $T$ is asymptotically contractive.
Then $T$ has a fixed point.
\end{thm}

\section{Examples}
\label{SC:examples}

In Proposition \ref{PROP:conditions},
 we prove
 \ref{ENUM:CONDI:asympt-cont-a}
 $\Rightarrow$ \ref{ENUM:CONDI:infinity-a}
 $\Rightarrow$ \ref{ENUM:CONDI:for-afpp}
 $\Rightarrow$ \ref{ENUM:CONDI:basic}.
In this section,
 we give three examples
 which show that the inverse of the above implications
 do not hold in general.

\begin{exmp}
\label{EX:ac-infi}
Put $E = \RealNumber$ and $C = [1,\infty)$.
Define a nonexpansive mapping $T$ on $C$ by
 $$ Tx = x - \log(x) $$
 for all $x \in C$.
Then $T$ satisfies \ref{ENUM:CONDI:infinity-a}
 and does not satisfy \ref{ENUM:CONDI:asympt-cont-a}.
\end{exmp}

\begin{proof}
Since
 $$ \toMUGEN{\limsup}{y}{C}
  \frac{\| T 1 - T y \|}{\| 1 - y \|}
 =
 \lim_{y \rightarrow \infty}
  \frac{y - \log(y) - 1}{y - 1}
 =
 1 , $$
 $T$ does not satisfy \ref{ENUM:CONDI:asympt-cont-a}
 by Proposition \ref{PROP:unbounded}.
Since
 \begin{align*}
 \toMUGEN{\lim}{y}{C}
  \Big( \| T 1 - T y \| - \| 1 - y \| \Big)
 &=
 \lim_{y \rightarrow \infty}
  \Big( \big( y - \log(y) - 1 \big) - \big( y - 1 \big) \Big) \\*
 &=
 \lim_{y \rightarrow \infty}
  - \log(y)
 =
 - \infty ,
 \end{align*}
 $T$ satisfies \ref{ENUM:CONDI:infinity-a}
 by Proposition \ref{PROP:unbounded}.
\end{proof}

\begin{exmp}
\label{EX:infi-afpp}
Let $E = c_0$ be the Banach space
 consisting of all real sequences converging to $0$
 with supremum norm.
Define a closed convex subset $C$ of $E$ by
 $$ C =
 \{ x \in E : 0 \leq x(n) \leq n \text{ for all } n \in \NaturalNumber \} .$$
Define a nonexpansive mapping $T$ on $C$ by
 $$ (Tx)(n) = \max\{ 0, x(n)-2 \} $$
 for $n \in \NaturalNumber$.
Then $T$ satisfies \ref{ENUM:CONDI:for-afpp}
 and does not satisfy \ref{ENUM:CONDI:infinity-a}.
\end{exmp}

\begin{proof}
Put $\lambda = 3$ and $x_1 = 0 \in C$.
It is clear that $T x_1 = 0$.
Fix $\eta > 0$ and choose $n \in \NaturalNumber$
 with $\eta < n$ and $2 \leq n$.
Put $y \in C$ by
 $$ y(k)
 =
 \begin{cases}
 n, & \text{if } k = n, \\
 0, & \text{if } k \neq n .
 \end{cases} $$
Then
 $ \| y \| = n > \eta $ and
 $$ (Ty) (k)
 =
 \begin{cases}
 n-2, & \text{if } k = n, \\
 0, & \text{if } k \neq n .
 \end{cases} $$
So, we have
 $$  \| T x_1 - T y \|
 = \| T y \| = n - 2
 > n - \lambda
 = \| y \| - \lambda
 = \| x_1 - y \| - \lambda . $$
Therefore $T$ does not satisfy \ref{ENUM:CONDI:infinity-a}.
We next put $\lambda = 1$ and fix $x_1 \in C$.
Then there exists $n_1 \in \NaturalNumber$ such that
 $0 \leq x_1(n) < 1$ for all $n \in \NaturalNumber$ with $n \geq n_1$.
By the definition of $T$,
 $(Tx_1)(n) = 0$ for $n \in \NaturalNumber$ with $n \geq n_1$.
Put $\eta = n_1 + 5$, and
 fix $y \in C$ with $\| y \| > \eta$.
We choose $n_2 \in N$ with $y(n_2) = \| y \|$.
Then from the definition of $C$, we have
 $$ n_1 < n_1 + 5 = \eta <  \| y \| = y(n_2) \leq n_2 . $$
It is clear that $\| y \| = y(n_2) > 2$.
For $n \in \NaturalNumber$ with $n < n_1$, we have
 $$ | \; (T x_1)(n) - (Ty)(n) \; | \leq n < n_1 < n_1 + 3 = \eta - 2 . $$
On the other hand,
 for $n \in \NaturalNumber$ with $n \geq n_1$, we have
 $$ | \; (T x_1)(n) - (Ty)(n) \; | = | \; (Ty)(n) \; | = \max\{ y(n)-2, 0 \}
 \leq \| y \| - 2 = y(n_2) - 2 . $$
Since $n_1 < n_2$, $\eta - 2 < \| y \| - 2 = y(n_2)-2$, and
 $$ | \; (T x_1)(n_2) - (Ty)(n_2) \; | = \max\{ y(n_2)-2, 0 \}
 = y(n_2)-2 , $$
 we have
 $$ \| T x_1 - Ty \| = y(n_2)-2 . $$
So, we obtain
 \begin{align*}
 \| T x_1 - T y \|
 &= y(n_2) - 2 \\*
 &\leq y(n_2) - x_1(n_2) - \lambda \\*
 &\leq \| x_1 - y \| - \lambda .
 \end{align*}
This implies \ref{ENUM:CONDI:for-afpp}.
This completes the proof.
\end{proof}

\begin{exmp}
\label{EX:afpp-basic}
Put $E = \RealNumber$ and $C = [1,\infty)$.
Define a nonexpansive mapping $T$ on $C$ by
 $$ Tx = x $$
 for all $x \in C$.
Then $T$ satisfies \ref{ENUM:CONDI:basic}
 and does not satisfy \ref{ENUM:CONDI:for-afpp}.
\end{exmp}

\begin{proof}
Since
 $$ \| T x - T y \| = \| x - y \| = \| x - y \| - \| x - T x\| $$
 for all $x, y \in C$,
 $T$ satisfies \ref{ENUM:CONDI:basic}.
And from the first equality,
 $T$ does not satisfy \ref{ENUM:CONDI:for-afpp}.
\end{proof}

\end{document}